\documentclass{article}
\usepackage{stmaryrd}
\usepackage{mathrsfs}
\usepackage[centertags]{amsmath}
\usepackage{amsfonts, dsfont}
\usepackage{amssymb}
\usepackage{amsthm}
\usepackage{color}

\input xy
\xyoption{all}

\theoremstyle{plain}
\newtheorem{thm}{Theorem}[section]
\newtheorem{cor}[thm]{Corollary}
\newtheorem{lem}[thm]{Lemma}
\newtheorem{prop}[thm]{Proposition}

\theoremstyle{definition}
\newtheorem{defn}[thm]{Definition}
\newtheorem{remark}[thm]{Remark}
\newtheorem*{ack}{Acknowledgments}

\newcommand{\bd}{\begin{defn}}
\newcommand{\ed}{\end{defn}}
\newcommand{\bl}{\begin{lem}}
\newcommand{\el}{\end{lem}}
\newcommand{\bp}{\begin{prop}}
\newcommand{\ep}{\end{prop}}
\newcommand{\bt}{\begin{thm}}
\newcommand{\et}{\end{thm}}
\newcommand{\bc}{\begin{cor}}
\newcommand{\ec}{\end{cor}}
\newcommand{\br}{\begin{remark}}
\newcommand{\er}{\end{remark}}
\newcommand{\bdi}{\begin{diagram}}
\newcommand{\edi}{\end{diagram}}
\newcommand{\beq}{\begin{equation}}
\newcommand{\eeq}{\end{equation}}
\newcommand{\ba}{\begin{array}}
\newcommand{\ea}{\end{array}}
\newcommand{\bpf}{\begin{proof}}
\newcommand{\epf}{\end{proof}}

\newcommand{\Z}{\mathds{Z}}
\newcommand{\Q}{\mathds{Q}}
\newcommand{\Zp}{\mathds{Z}_{p}}
\newcommand{\Qp}{\mathds{Q}_{p}}
\newcommand{\al}{\alpha}

\newcommand{\Ga}{\Gamma}

\newcommand{\la}{\lambda}

\DeclareMathOperator{\Sel}{Sel} \DeclareMathOperator{\Gal}{Gal}
\DeclareMathOperator{\Hom}{Hom} \DeclareMathOperator{\rank}{rank}
\DeclareMathOperator{\Ext}{Ext} \DeclareMathOperator{\Ind}{Ind}

\newcommand{\cyc}{\mathrm{cyc}}
\newcommand{\cts}{\mathrm{cts}}
\newcommand{\Hi}{H_{\mathrm{Iw}}}

\newcommand{\Op}{\mathcal{O}}
\newcommand{\m}{\mathfrak{m}}

\newcommand{\ot}{\otimes}
\newcommand{\ilim}{\displaystyle \mathop{\varinjlim}\limits}
\newcommand{\plim}{\displaystyle \mathop{\varprojlim}\limits}
\newcommand{\im}{\mathrm{im}\,}

\newcommand{\lra}{\longrightarrow}

\newcommand{\ps}[1]{\llbracket #1 \rrbracket}

\numberwithin{equation}{section}

\begin{document}
\title{Fine Selmer groups of congruent Galois representations}
\author{Meng Fai Lim\footnote{School of Mathematics and Statistics $\&$ Hubei Key Laboratory of Mathematical Sciences,
Central China Normal University, Wuhan, 430079, P.R.China.
 E-mail: \texttt{limmf@mail.ccnu.edu.cn}} \quad
 Ramdorai Sujatha\footnote{Department of Mathematics,
 University of British Columbia,
 Vancouver BC, V6T 1Z2,
 Canada.
 E-mail: \texttt{sujatha@math.ubc.ca}}}
\date{}
\maketitle

\begin{abstract} \footnotesize
\noindent In this paper, we study the fine Selmer groups of two
congruent Galois representations over an admissible $p$-adic Lie
extension. We show that under appropriate congruence conditions, if
the dual fine Selmer group of one is pseudo-null, so is the other.
Our results also compare the $\pi$-primary submodules of the two
dual fine Selmer groups. We then apply our results to compare the
structure of Galois group of the maximal abelian unramified pro-$p$
extension of an admissible $p$-adic Lie extension and the structure
of the dual fine Selmer group over the said admissible $p$-adic Lie
extension.


\end{abstract}

\section{Introduction}
Throughout the paper, $p$ will always denote a fixed prime. Let $F$
be a number field. If $p=2$, we assume that $F$ is totally
imaginary. We write $F^{\cyc}$ for the cyclotomic $\Zp$-extension of
$F$, whose Galois group $\Gal(F^{\cyc}/F)$ is in turn denoted $\Ga$.
Denote by $K(F^{\cyc})$ the maximal unramified pro-$p$ extension of
$F^{\cyc}$ in which every prime of $F^{\cyc}$ above $p$ splits
completely. In \cite{Iw2}, Iwasawa proved  that
$\Gal(K(F^{\cyc})/F^{\cyc})$ is a finitely generated torsion
$\Zp\llbracket\Ga\rrbracket$-module, and he further conjectured that
this Galois group is finitely generated over $\Zp$ (see also
\cite{Iw}). Throughout this article, we shall call this conjecture
the \textit{Iwasawa $\mu$-conjecture}. On the other hand, it has
been known that this finite generation property \textit{does not}
hold for the dual of the classical Selmer group of an abelian
variety over the cyclotomic $\Zp$-extension in general (see \cite[\S
10, Example 2]{Mazur}). It was only about a decade ago that Coates
and the second named author \cite{CS,Su} gave a correct formulation
of the analogue of the Iwasawa $\mu$-conjecture for an elliptic
curve. Namely, they considered a smaller group, called the fine
Selmer group, which is a subgroup of the classical Selmer group, and
they conjectured that the Pontryagin dual of this fine Selmer group
over $F^{\cyc}$ is finitely generated over $\Zp$ \cite[Conjecture
A]{CS}. Since then, analogues of this conjecture have been
formulated for fine Selmer groups attached to more general Galois
representations (see \cite{A,Jh,JhS,LimFine,LimMurty}). In this
paper, we shall collectively (and loosely) address these conjectures
as \textit{Conjecture A}. A striking observation is that, besides
being a natural analogue of the Iwasawa $\mu$-conjecture,
Conjecture A is related to the latter conjecture in a very
precise manner (see \cite[Theorem 3.4]{CS}, \cite[Theorem
3.5]{LimFine}, \cite[Theorem 5.5]{LimMurty}, \cite[Theorem 4.5]{Su}
or \cite[Section 8]{Wu}; also see Theorem \ref{fg cong corollary}
below).

In their paper \cite{CS}, Coates and the second author also studied
the structure of the fine Selmer group over extensions of $F$ whose
Galois group $G = \Gal(F_\infty/F)$ is a $p$-adic Lie group of
dimension larger than $1$. There they formulated an important
conjecture on the structure of the Pontryagin dual of the fine
Selmer group of an abelian variety which predicts that the said
module is pseudo-null over the Iwasawa algebra $\Zp\ps{G}$ (see
\cite[Conjecture B]{CS}). To some extent, their conjecture can be
thought as an analogue of a conjecture of Greenberg, which we now
briefly describe. Recall that a Galois extension $F_{\infty}$ of $F$
is said to be a strongly admissible, pro-$p$, $p$-adic Lie extension
of $F$ if (i) $G=\Gal(F_{\infty}/F)$ is  a compact pro-$p$, $p$-adic Lie
group without $p$-torsion, (ii) $F_{\infty}$ contains the cyclotomic
$\Zp$ extension $F^{\cyc}$ of $F$ and (iii) $F_{\infty}$ is
unramified outside a finite set of primes. We denote by $H$ the Galois group
$\Gal(F_{\infty}/F^{\cyc}).$ Let $K(F_{\infty})$
denote the maximal unramified abelian pro-$p$ extension of
$F_{\infty}$ in which every prime above $p$ splits completely. When
$F_{\infty}$ is the composite of all the $\Zp$-extensions of $F$,
Greenberg \cite{Gr} conjectured that
$\Gal(K(F_{\infty})/F_{\infty})$ is pseudo-null over $\Zp\ps{G}$.
(Actually, to be more precise, Greenberg's original conjecture is
concerned with the pseudo-nullity of a slightly bigger Galois
group.) For a general $F_{\infty}$, the validity of the
pseudo-nullity of $\Gal(K(F_{\infty})/F_{\infty})$ is not
guaranteed, and this was first observed by Hachimori and Sharifi in
\cite{HS}, where they constructed a class of  extensions  $F_{\infty}$ whose
Galois group $\Gal(K(F_{\infty})/F_{\infty})$ is not pseudo-null.
Despite these constructions of Hachimori and Sharifi, Coates and the
second named author have expressed optimism that the corresponding
assertion for the dual fine Selmer group of an  elliptic curve
should hold, for strongly admissible extensions
$F_{\infty}$ (see \cite[Section 4]{CS}).

Since then, the question of the pseudo-nullity of the dual fine
Selmer group  over a general, strongly admissible  $p$-adic
Lie extension has been the subject of much study (see
\cite{Bh,Jh,LimPS,LimFine,Oc}). In \cite{Jh}, Jha formulated an
analogue of this conjecture for a Hida family and its
specialization. In his work, Jha was able to show that if the dual
fine Selmer group of one specialization of the Hida family is
pseudo-null, then the dual fine Selmer groups for all but finitely
many specializations of the Hida family are also pseudo-null (see
\cite[Theorem 10]{Jh}). Of course, in view of the conjectures of
Coates-Sujatha and Jha, one expects that every (arithmetic)
specialization has a pseudo-null dual fine Selmer group, and
therefore, the theorem of Jha gives a very strong evidence to this.
A careful examination of Jha's proof actually yields a criterion of
determining which arithmetic specialization has a pseudo-null dual
fine Selmer group, and the determining criterion relies on the
structure of the central torsion submodule of the dual fine Selmer
group of the big Galois representation. In view of Jha's result, it
is then natural to investigate the preservation of pseudo-nullity of
the fine Selmer groups of congruent Galois representations in
general. More precisely, motivated by Jha's theorem, the following
natural question  is  of interest: \textit{Suppose that the dual fine
Selmer group of one of the congruent Galois representations is
pseudo-null, can one deduce the pseudo-nullity of the fine Selmer
group of the other representation via a criterion on the structure
of the dual fine Selmer group of the initial representation?} The
primary goal of this paper is to develop such a criterion. As will
be seen below, our criterion depends on the
$p$-primary submodule of the dual fine Selmer group of the initial
representation (see Theorems \ref{ps cong} and \ref{ps cong high}).
We emphasise that these dual fine Selmer groups are expected to have
trivial $\mu_{\Op\ps{G}}$-invariants and hence their $p$-primary
submodules are pseudo-null as $\Op\ps{G}$-modules. However,
as $\Op\ps{H}$-modules, the structures of these $p$-primary
submodules are not known and it is not clear what to expect of them.
Our results (Theorems \ref{ps cong} and \ref{ps cong high}) will
therefore consist of considering the situations when the dual fine
Selmer group of the initial representation has a trivial
$\mu_{\Op\ps{H}}$-invariant and when it does not.

We then apply our criterion to study the relation between
$\Gal(K(F_{\infty})/F_{\infty})$ and the dual fine Selmer groups.
Motivated by the relation between the Iwasawa $\mu$-conjecture and
Conjecture A, one may ask whether there is an analogous relationship
between the pseudo-nullity of $\Gal(K(F_{\infty})/F_{\infty})$ and
the pseudo-nullity of the dual fine Selmer groups. Of course, the
constructions of Hachimori-Sharifi tell us that such an analogue does not
hold on the nose. Nevertheless, we can still ask the question of deducing the
pseudo-nullity of the dual fine Selmer groups from the knowledge of
the pseudo-nullity of the Galois group of the maximal abelian
unramified pro-$p$ extension (see Question B$'$ in Section
\ref{Galois and fine}). Some partial result in this direction has
been obtained by the first named author in \cite[Theorem
2.3]{LimPS}. In this paper, we give a refinement of these results (see
Propositions \ref{pseudo-null main} and \ref{ps cong high Galois}).
We also relate the $p$-primary submodule of
$\Gal(K(F_{\infty})/F_{\infty})$ and the $\pi$-primary submodule of
the dual fine Selmer group (as $\Zp\ps{H}$-modules).

Finally, we return to the situation of an elliptic curve $E$ with good ordinary reduction at all primes above $p$. We deduce a relation between Conjecture A and the structure
of the Selmer group of the said elliptic curve (see Theorem \ref{Sel vs fine}). This relation can be thought of as the mod-$p$ analogue of Mazur's conjecture, which states that the dual Selmer group over the cyclotomic $\Zp$-extension is a torsion module over the Iwasawa algebra $\Zp\ps{\Ga}$. It is well known that this statement is equivalent to the defining sequence for the Selmer group being short exact and the validity of an appropriate version of the Weak Leopoldt conjecture. In fact, the corresponding equivalence between the dual Selmer group being torsion over the associated Iwasawa algebra and the defining sequence for the Selmer group being short exact (modulo an appropriate version of the Weak Leopoldt conjecture) is true for more general strongly admissible $p$-adic Lie extensions (see \cite[Theorem 4.12]{C99} or \cite[Proposition 3.3]{LimMHG}). In our situation, Theorem \ref{Sel vs fine} asserts that the dual strict Selmer group of $E[p]$ being torsion over the associated Iwasawa algebra over $\mathbb{F}_p$ is equivalent to the defining sequence for the strict Selmer group of $E[p]$ being short exact and the validity of Conjecture A. (Here the strict Selmer group is in the sense of Greenberg \cite{G89}.) This result provides conceptual clarity to a result of
Vatsal and Greenberg (see \cite[Page 18, Statement A]{GV}). In particular, we shall apply Theorem \ref{Sel vs fine} to give a proof of the said result of Vatsal and Greenberg (see Corollary \ref{GV proof}).

We now give a brief description of the layout of the paper. In
Section \ref{Algebraic Preliminaries}, we recall certain algebraic
notions which will be used subsequently in the paper. In Section
\ref{fine Selmer group}, we introduce the fine Selmer groups. It is
also here that we establish our main results. In Section \ref{Galois
and fine}, we then apply the criterion developed in Section
\ref{fine Selmer group} to relate the Iwasawa module theoretical
structure of Galois groups and fine Selmer groups. In Section \ref{Galois
and fine}, we then apply the criterion developed in Section
\ref{fine Selmer group} to relate the Iwasawa module theoretical
structure of Galois groups and fine Selmer groups. In Section \ref{further remark section}, we revisit
the situation of an elliptic curve.

\begin{ack}
Some part of the research of this article took place when the first
author was visiting the National Center for Theoretical Sciences,
the Institute of Mathematics of Academia Sinica, the University of
Toronto and the Ganita Lab. The first author would like to thank
these institutes for their hospitality. The authors also
thank Manfred Kolster and Romyar Sharifi for many interesting
conversations on the subject of $p$-primary submodules of Iwasawa
modules. Finally, the first author's research is supported by the
National Natural Science Foundation of China under The Research Fund
for International Young Scientists (Grant No: 11550110172), and the
second author gratefully acknowledges support from  NSERC
(Discovery Grant No: 402071).

 \end{ack}

\section{Algebraic Preliminaries} \label{Algebraic Preliminaries}

As before, $p$ will denote a prime number.
Let $\Op$ be the ring of integers of a fixed finite extension of $\Qp$.
For a compact $p$-adic Lie group $G$, the
completed group algebra of $G$ over $\Op$ is defined by
 \[ \Op\ps{G} = \plim_U \Op[G/U], \]
where $U$ runs over the open normal subgroups of $G$ and the inverse
limit is taken with respect to the canonical projection maps.
Throughout the paper, we usually work under the assumption that our group $G$
is pro-$p$ and has no $p$-torsion. In this setting, it
is known that $\Op\ps{G}$ is an Auslander regular ring with no zero divisors (cf.
\cite[Theorem 3.26]{V02} and \cite{Neu}). Therefore, the ring $\Op\ps{G}$
admits a skew field $Q(G)$ which is flat over $\Op\ps{G}$ (see
\cite[Chapters 6 and 10]{GW} or \cite[Chapter 4, \S 9 and \S
10]{Lam}). The $\Op\ps{G}$-rank of a finitely generated
$\Op\ps{G}$-module $M$ is then defined to be
$$ \rank_{\Op\ps{G}}(M)  = \dim_{Q(G)} \big(Q(G)\ot_{\Op\ps{G}}M\big). $$
We say that the $\Op\ps{G}$-module $M$ is \textit{torsion} if
$\rank_{\Op\ps{G}} M = 0$. We frequently make use of the following
well-known equivalent characterization for a torsion
$\Op\ps{G}$-module without further comment, namely:
$\Hom_{\Op\ps{G}}(M, \Op\ps{G})=0$ (for instance, see \cite[Lemma
4.2]{LimFine}). A finitely generated torsion $\Op\ps{G}$-module $M$
is then said to be \textit{pseudo-null} if $\Ext^1_{\Op\ps{G}}(M,
\Op\ps{G}) =0$.

Let $\pi$ denote a fixed local parameter for $\Op$ and let $k$
denote the residue field of $\Op$. The completed group algebra of
$G$ over $k$ is denoted $k\ps{G}$ and we have, as before,
 \[ k\ps{G} = \plim_U k[G/U], \]
where $U$ runs over the open normal subgroups of $G$ and the inverse
limit is taken with respect to the canonical projection maps. For a
compact $p$-adic Lie group $G$ without $p$-torsion, it follows from
\cite[Theorem 3.30(ii)]{V02} (or \cite[Theorem A.1]{LimFine}) that
$k\ps{G}$ is an Auslander regular ring. Furthermore, if $G$ is
pro-$p$, then the ring $k\ps{G}$ has no zero divisors (cf.
\cite[Theorem C]{AB}). Therefore, one can define the notion of
$k\ps{G}$-rank as above when $G$ is pro-$p$ without $p$-torsion.
Similarly, we say that a finitely generated  $k\ps{G}$-module $N$ is a
\textit{torsion} module if $\rank_{k\ps{G}}N = 0$.

For a given finitely generated $\Op\ps{G}$-module $M$, we denote by
$M(\pi)$ the $\Op\ps{G}$-submodule of $M$ which consists of elements
of $M$ that are annihilated by some power of $\pi$. We shall call
$M(\pi)$ the $\pi$-primary submodule of $M$. A finitely generated
$\Op\ps{G}$-module is then said to be $\pi$-primary if $M=M(\pi)$.
Since the ring $\Op\ps{G}$ is Noetherian, the $\pi$-primary
submodule $M(\pi)$ of $M$ is also finitely generated over
$\Op\ps{G}$. Therefore, one can find an integer $r\geq 0$ such that
$\pi^r$ annihilates $M(\pi)$. Following \cite[Formula (33)]{Ho}, we
define the \textit{$\mu$-invariant}
  \[\mu_{\Op\ps{G}}(M) = \sum_{i\geq 0}\rank_{k\ps{G}}\big(\pi^i
   M(\pi)/\pi^{i+1}\big). \]
(For another alternative, but equivalent, definition, see
\cite[Definition 3.32]{V02}.) By the above discussion and our
definition of $k\ps{G}$-rank, the sum on the right is a finite one.
It is immediate from the definition that $\mu_{\Op\ps{G}}(M) =
\mu_{\Op\ps{G}}(M(\pi))$.

Now suppose that $G$ is pro-$p$ without $p$-torsion. Then as seen in
the discussion above, both rings $\Op\ps{G}$ and $k\ps{G}$ are
Auslander regular and have no zero divisors. For a finitely
generated $\Op\ps{G}$-module $M$, it then follows from either
\cite[Proposition 1.11]{Ho2} or \cite[Theorem 3.40]{V02} that there
is a $\Op\ps{G}$-homomorphism
\[ \varphi: M(\pi) \lra \bigoplus_{i=1}^s\Op\ps{G}/\pi^{\al_i},\] whose
kernel and cokernel are pseudo-null $\Op\ps{G}$-modules, and where
the integers $s$ and $\al_i$ are uniquely determined. We call
$\bigoplus_{i=1}^s\Op\ps{G}/\pi^{\al_i}$  the \textit{elementary
representation} of $M(\pi)$. In fact, in the process of establishing
the above, one also has that $\mu_{\Op\ps{G}}(M) = \displaystyle
\sum_{i=1}^s\al_i$ (see loc. cit.). We will set
\[\theta_{\Op\ps{G}}(M) : = \max_{1\leq i\leq s}\{\al_i\}.\]

For convenience, we recall some of the necessary results from
\cite{LimCMu} that we will need here.

 \bp \label{pseudo-isomorphic}
 Let $M$ and $N$ be two finitely generated $\Op\ps{G}$-modules such
that $M$ is a torsion $\Op\ps{G}$-module and such that
$\mu_{\Op\ps{G}}(M/\pi^i) = \mu_{\Op\ps{G}}(N/\pi^i)$ for every
$1\leq i \leq \theta_{\Op\ps{G}}(M)+1$.

Then $N$ is torsion over $\Op\ps{G}$ and we have the equality
$\theta_{\Op\ps{G}}(M) = \theta_{\Op\ps{G}}(N)$. Furthermore,
$M(\pi)$ and $N(\pi)$ have the same elementary representations as $\Op\ps{G}$-modules.
 \ep

We actually require the following pseudo-null analogue of the above proposition.

\bp \label{pseudo-isomorphic2} Let $H$ be a closed normal subgroup
of $G$ with $G/H\cong\Zp$. Let $M$ and $N$ be two finitely generated
$\Op\ps{G}$-modules which are also finitely generated over
$\Op\ps{H}$. Suppose that $M$ is a pseudo-null $\Op\ps{G}$-module
and that $\mu_{\Op\ps{H}}(M/\pi^i) = \mu_{\Op\ps{H}}(N/\pi^i)$
for every $1\leq i \leq \theta_{\Op\ps{H}}(M)+1$.

Then $N$ is pseudo-null over $\Op\ps{G}$ and we have the equality
$\theta_{\Op\ps{H}}(M) = \theta_{\Op\ps{H}}(N)$. Furthermore,
$M(\pi)$ and $N(\pi)$ have the same elementary representations as $\Op\ps{H}$-modules.
 \ep

\bpf
 If $M$ is an $\Op\ps{G}$-module which is finitely generated over
$\Op\ps{H}$, it then follows from a well-known result of Venjakob \cite{V03}
that $M$ is a pseudo-null $\Op\ps{G}$-module if and
only if it is a torsion $\Op\ps{H}$-module. The conclusion of the proposition is now immediate consequence
of Proposition \ref{pseudo-isomorphic}.
\epf

\br
Note that we are concerned with the $\mu_{\Op\ps{H}}$-invariants in the proposition which may not be zero even if the modules
in question are pseudo-null over $\Op\ps{G}$.
\er

We end this section with two more useful lemmas. For a closed
subgroup $U$ of a compact $p$-adic Lie group $H$ and a compact
$\Op\ps{U}$-module $M$, we denote by
% Changes made on  4/2/2017
$\Ind^H_U(M) :=
M\hat{\ot}_{\Op\ps{U}}\Op\ps{H}$ the compact induction of $M$ from
$U$ to $H$. If $M$ is a compact $k\ps{U}$-module, we can also view
it as a compact $\Op\ps{U}$-module and we have the identification
$\Ind^H_U(M) = M\hat{\ot}_{k\ps{U}}k\ps{H}$.

\bl \label{ind} Let $H$ be a compact pro-$p$ $p$-adic Lie group
without $p$-torsion and $U$ a closed subgroup of $H$ of dimension at
least 1. Let $N$ be a finitely generated $\Op\ps{U}$-module which is
finitely generated over $\Op$, and hence torsion as an $\Op[[U]]$-module.
Then $\Ind^H_U(N)$ is a finitely
generated torsion $\Op\ps{H}$-module with
$\mu_{\Op\ps{H}}(\Ind^H_U(N)) = 0$. \el

\bpf
 By \cite[Lemma 5.5]{OcV02}, we have an isomorphism
 \[ \Hom_{\Op\ps{H}}(\Ind^H_U(N), \Op\ps{H})\cong
 \Ind^H_U\Big(\Hom_{\Op\ps{U}}(N, \Op\ps{U})\Big).\]
Since $U$ has dimension at least 1 and $N$ is finitely generated
over $\Op$, it follows that $N$ is torsion over $\Op\ps{U}$. Thus,
the term on the right of the isomorphism is zero, and by the
isomorphism, this in turns implies that $\Ind^H_U(N)$ is a finitely
generated torsion $\Op\ps{H}$-module. Again, using the fact that
$U$ has dimension at least 1 and that $N$ is finitely generated over
$\Op$, we have that $\pi^n N(\pi)/\pi^{n+1}$ is torsion over
$k\ps{U}$. Since $\Ind^H_U(-)$ is an exact functor (cf.\ \cite[Lemma
6.10.8]{RZ}), it follows from a straightforward argument that
\[\pi^n \left((\Ind^H_U
N)(\pi)\right)/\pi^{n+1} \cong \Ind^H_U \big(\pi^n
N(\pi)/\pi^{n+1}\big).\] Via a similar argument as above, we have
that $\pi^n \left((\Ind^H_U N)(\pi)\right)/\pi^{n+1}$ is a torsion
$k\ps{H}$-module. Combining this observation with the definition of
the $\mu_{\Op\ps{H}}$-invariant will yield that
$\mu_{\Op\ps{H}}(\Ind^H_U(N)) = 0$.
 \epf

The next lemma is a refinement of \cite[Lemma 2.4.3]{LimCMu}.

\bl \label{pseudo-isomorphism lemma}
 Suppose that $H$ is a compact pro-$p$ $p$-adic Lie group
without $p$-torsion. Let $\varphi: M \lra N$ be a homomorphism of
finitely generated $\Op\ps{H}$-modules, whose kernel and cokernel
are finitely generated torsion $\Op\ps{H}$-modules with trivial
$\mu_{\Op\ps{H}}$-invariants. Then $M(\pi)$ and $N(\pi)$ have the
same elementary representations as $\Op\ps{H}$-modules. \el

\bpf
 The statement will follow if it holds in the two special cases of exact
sequences
\[\ba{c} 0\lra P \lra M \lra N\lra 0, \\
  0\lra M \lra N \lra P\lra 0, \ea
\] where $P$ is a finitely generated $\Op\ps{H}$-module with trivial
$\mu_{\Op\ps{H}}$-invariant. We will prove the second case, the
first case has a similar argument. Choose a sufficiently large $n$
such that $\pi^n$ annihilates $M(\pi)$ and $N(\pi)$. Consider the
following commutative diagram
  \[ \SelectTips{eu}{}
\xymatrix{
  0 \ar[r] & M \ar[d]^{\pi^n} \ar[r]^{} & N \ar[d]^{\pi^n}
  \ar[r]^{} & P \ar[d]^{\pi^n} \ar[r]& 0 \\
  0 \ar[r] & M \ar[r] & N
  \ar[r] & P \ar[r] &0   }
\] with exact rows, and the vertical maps are given by
multiplication by $\pi^n$. From this, we obtain an exact sequence
\[ 0\lra M(\pi) \stackrel{\varphi}{\lra} N(\pi) \lra P(\pi). \]
Since $\mu_{\Op\ps{H}}(P)=0$, it follows from \cite[Remark
3.33]{V02} that $P(\pi)$ is a pseudo-null $\Op\ps{H}$-module. This
gives the required conclusion in the lemma.
% Changes made on  4/2/2017
Let
\[ f: N(\pi) \lra \bigoplus_{i=1}^s\Op\ps{G}/\pi^{\al_i}\]
 be a homomorphism of $\Op\ps{H}$-modules,
whose kernel and cokernel are pseudo-null $\Op\ps{H}$-modules. Then
 \[f\circ \varphi: M(\pi) \lra
\bigoplus_{i=1}^s\Op\ps{G}/\pi^{\al_i}\] is a homomorphism of
$\Op\ps{H}$-modules, whose kernel and cokernel are pseudo-null.
Therefore, $M(\pi)$ and $N(\pi)$ have the same elementary
representations.\epf

\section{Fine Selmer groups} \label{fine Selmer group}

We now move to arithmetic. As before, $p$ denotes a prime. Let
$F$ be a number field. If $p=2$, assume further that the number field $F$ has
no real primes. Let $K$ be a fixed finite
extension of $\Qp$, whose ring of integers is denoted  $\Op$. We shall also fix a choice of local parameter $\pi$ for $\Op$. Suppose that we are given
a finite dimensional $K$-vector
space $V$ with a continuous $\Gal(\bar{F}/F)$-action which is unramified
outside a finite set of primes. By a standard compactness argument, $V$ contains a $\Gal(\bar{F}/F)$-stable
$\Op$-lattice which we will fix once and for all,  and denote it by $T$.
Write $A= V/T$. Let $S$ denote a finite set of primes of $F$
which contains all the primes above $p$, the ramified primes of $A$ and the
infinite primes. We then denote by $F_S$ the maximal algebraic extension of
$F$ unramified outside $S$. For any algebraic (possibly infinite)
extension $\mathcal{L}$ of $F$ contained in $F_S$, we shall write
$G_S(\mathcal{L}) = \Gal(F_S/\mathcal{L})$.

Let $v$ be a prime in $S$. For each finite extension $L$ of $F$
contained in $F_S$, we define
 \[K_v^i(A/L) = \bigoplus_{w|v}H^i(L_w, A) \quad(i=0,1),\]
where $w$ runs over the (finite) set of primes of $L$ above $v$. If
$\mathcal{L}$ is an infinite extension of $F$ contained in $F_S$, we
define
\[K_v^i(A/\mathcal{L}) = \ilim_L K_v^i(A/L),\]
where the direct limit is taken over all finite extensions $L$ of
$F$ contained in $\mathcal{L}$ under the restriction maps on the cohomology.

For any algebraic (possibly infinite) extension $\mathcal{L}$ of $F$
contained in $F_S$, the \textit{fine Selmer group} of $A$ over
$\mathcal{L}$ (with respect to $S$) is defined to be
\[ R_S(A/\mathcal{L}) = \ker\Big(H^1(G_S(\mathcal{L}), A)\lra \bigoplus_{v\in S}K_v^1(A/\mathcal{L})
\Big). \] We shall write $Y_S(A/\mathcal{L})$ for the Pontryagin
dual $R_S(A/\mathcal{L})^{\vee}$ of the fine Selmer group. The
following conjecture, first formally stated in \cite{CS}, which is also implicit in \cite{Iw,Iw2} (see discussion below), is now folklore.

\medskip \noindent \textbf{Conjecture A.} \textit{For
any number field $F$,\, $Y_S(A/F^{\cyc})$ is a finitely generated
$\Op$-module.}

\medskip
 We take this opportunity to highlight  the following observation. For any
extension $\mathcal{L}$ of $F^{\cyc}$ contained in $F_S$, let
$K(\mathcal{L})$  be the maximal unramified pro-$p$ extension of
$\mathcal{L}$ where every prime of $\mathcal{L}$ above $p$ splits
completely. It then follows that every
finite prime of $\mathcal{L}$ splits completely in
$K(\mathcal{L})$. Therefore, in the case when $V = \Qp$ and $A=\Qp/\Zp$, the dual
of the fine Selmer group $Y_S(A/F^{\cyc})$ is precisely
$\Gal(K(F^{\cyc})/F^{\cyc})$. In this context, Conjecture A (for $A=\Qp/\Zp$) is equivalent to
the conjecture made by Iwasawa \cite{Iw, Iw2} which we shall call
the Iwasawa $\mu$-conjecture for $F^{\cyc}$.

For the purpose of the paper, we will require another equivalent
formulation of Conjecture A. Write $T^*=
\Hom_{\cts}(A,\mu_{p^{\infty}})$,
where $\mu_{p^{\infty}}$ denotes
the group of all $p$-power roots of unity. (One should not confuse
this with $T$.) In fact, it's not difficult to show that $T^* =
\Hom_{\Op}(T,\Op(1))$, where `(1)' means Tate twist. The relationship between $T, T^*$ and $A$ is illustrated by the following diagram
\[  \entrymodifiers={!! <0pt, 0.8ex>+} \SelectTips{eu}{}\xymatrix@R=.6in @C=1in{
     T \ar[d]_{-\ot_{\Op}K/\Op} \ar[r]^{\Hom_{\Op}(-, \Op(1))} &
    T^* \ar[l] \ar[dl]
     \\
    A  \ar[ur]_{\Hom_{\cts}(-,\mu_{p^{\infty}})}
    &  }  \]  By an
application of the Poitou-Tate duality (see \cite[Eqn. (45)]{CS}), we have the following exact
sequence
 \[
0\lra Y_S(A/\mathcal{L}) \lra \Hi^2(\mathcal{L}/F, T^*)\lra
\Bigg(\bigoplus_{v\in S} K_v^0(A/\mathcal{L})\Bigg)^{\vee},
\]
 where $\Hi^2(\mathcal{L}/F, T^*) =
\plim_{L}H^2(G_S(L), T^*)$. By \cite[Lemma 3.2]{CS} (or see
\cite[Lemma 3.4]{LimFine}), Conjecture A holds for $A$ over
$F^{\cyc}$ if and only if $\Hi^2(F^{\cyc}/F, T^*)$ is finitely
generated over $\Op$. We will frequently make use of this equivalent
formulation of Conjecture A without further comment.

As mentioned in the introductory section, Conjecture A is intimately
related to the Iwasawa $\mu$-conjecture. This was observed by
Coates and the second author in the context of an  elliptic
curve (see \cite[Theorem 3.4]{CS}). Subsequently, this observation
has been generalized to more general $p$-adic representations (for
instance, see \cite[Theorem 8]{JhS} and \cite[Theorem
3.5]{LimFine}). We now recall this result. Let $\bar{\rho} :
G_{S}(F)\lra \mathrm{Aut}_k(E[\pi])$ be the representation induced
by the Galois action of $G_S(F)$ on $E[\pi]$. Denote by $F(E[\pi])$
the subextension of $F_S$ which is fixed by the kernel of
$\bar{\rho}$. Note that this is a finite Galois extension of $F$.
The following is a slight refinement of the  observation of
Coates and the second author.

\bt \label{fg cong corollary}
 Let $F^{\cyc}$ be the cyclotomic $\Zp$-extension of $F$.
 Suppose that $F(E[\pi])$ is a finite $p$-extension of $F$.
  Then Conjecture A holds for $E$ over $F^{\cyc}$ if and only if
  the Iwasawa $\mu$-conjecture holds for $F^{\cyc}$.
\et

There are a number of approaches towards proving the above theorem
(see \cite[Theorem 3.4]{CS}, \cite[Theorem 8]{JhS}, \cite[Theorem
3.5]{LimFine}, \cite[Theorem 5.5]{LimMurty} and \cite[Theorem
4.5]{Su}). We will give a proof of the theorem along the lines of
\cite[Theorem 8]{JhS} and \cite[Theorem 3.5]{LimFine}. To prepare
for the proof, we need to introduce another Galois representation.
Let $V'$ be another finite dimensional $K$-vector space with a
continuous $\Gal(\bar{F}/F)$-action which is unramified outside a
finite set of primes. Write $B=V'/T'$ for a (fixed) Galois invariant
lattice $T'$ of $V'$. From now on, our (finite) set $S$ of primes is
always assumed to contain  the primes above $p$, the ramified primes
of $A$ and $B$, and the infinite primes. We can now state the
following.

\bp \label{fg cong} Suppose that there is an isomorphism $A[\pi]\cong B[\pi]$ of
$G_S(F)$-modules. Then Conjecture A holds for $Y_S(A/F^{\cyc})$ if
and only if Conjecture A holds for $Y_S(B/F^{\cyc})$. \ep

\bpf For $Z=A,B$, by a similar argument to that in  \cite[Proposition 4.6]{Su}, one can show that
 Conjecture A is equivalent to the assertion
that $H^2(G_S(F^{\cyc}), Z[\pi])=0$. The conclusion is now immediate from the congruence condition.
% Changes made on 4/2/2017
 \epf

We now give the proof of Theorem \ref{fg cong corollary}.

 \bpf[Proof of Theorem \ref{fg cong corollary}]
   Suppose, for now, that $L'$ is a finite $p$-extension of $L$. By a classical argument
(cf. \cite[Theorem 3]{Iw}), we have that the Iwasawa $\mu$-conjecture holds for $L^{\cyc}$
if and only if the Iwasawa $\mu$-conjecture holds for
$L'^{\cyc}$. Via a similar argument to that in \cite[Theorem 5.5]{LimMurty}, we also have that
 Conjecture A
holds for $Y_S(A/L^{\cyc})$ if and only if Conjecture A holds for
$Y_S(A/L'^{\cyc})$. Hence it suffices to prove the equivalence in
the theorem over $F(A[\pi])$. In particular, without loss of
generality, we may assume that $F=F(A[\pi])$. It then follows from
this that $A[\pi]$ is a trivial $G_S(F)$-module and that
$A[\pi]\cong (\Op/\pi)^d \cong (\Z/p)^{fd}$ as $G_S(F)$-modules.
Here $d=\mathrm{corank_{\Op}}(A)$ and $f=[K:\Qp]$. By Proposition
\ref{fg cong}, Conjecture A holds for $Y_S(A/F^{\cyc})$ if and only
if Conjecture A holds for $Y_S\left((\Qp/\Zp)^{fd}/F^{\cyc}\right)$.
But the latter clearly holds if and only if the Iwasawa
$\mu$-conjecture holds for $F^{\cyc}$. \epf

 We now turn to the situation of a higher dimensional
$p$-adic Lie extension. Recall that a Galois extension $F_{\infty}$ of $F$ is
said to be a strongly $S$-admissible, pro-$p$ , $p$-adic Lie extension of $F$ if
(i) $\Gal(F_{\infty}/F)$ is a compact  pro-$p$, $p$-adic Lie group
without $p$-torsion, (ii) $F_{\infty}$ contains $F^{\cyc}$ and (iii)
$F_{\infty}$ is contained in $F_S$, the maximal unramified outside $S$-extension of $F$,
where $S$ , as before, is a finite set of primes
which contains all the primes above $p$, the ramified primes fo rthe Galois representation,
and the infinite primes.. We write $G=\Gal(F_{\infty}/F)$ and
$H=\Gal(F_{\infty}/F^{\cyc})$. To facilitate further discussion, we
record the following lemma.

\bl \label{fg La H} Let $F_{\infty}$ be a strongly $S$-admissible
 pro-$p$, $p$-adic Lie extension of $F$. Then the following statements are
equivalent.

\begin{enumerate}
\item[$(a)$] $Y_S(A/F^{\mathrm{cyc}})$ is a finitely generated
$\Op$-module.

\item[$(b)$] $\Hi^2(F^{\cyc}/F, T^*)$ is a finitely generated $\Op$-module.

\item[$(c)$] $Y_S(A/F_{\infty})$ is a finitely generated
$\Op\ps{H}$-module.

\item[$(d)$] $\Hi^2(F_{\infty}/F, T^*)$ is a finitely generated $\Op\ps{H}$-module.
 \end{enumerate}\el

\bpf The proof  is entirely similar to that in \cite[Lemma 3.2]{CS}. \epf

By a well-known result of Venjakob \cite{V03},
we have that  an $\Op\ps{G}$-module $M$ which is finitely generated over
$\Op\ps{H}$,  is a pseudo-null $\Op\ps{G}$-module if and
only if it is a torsion $\Op\ps{H}$-module. In view of this, we may now pose the following question.

\medskip

\noindent \textbf{Question B:}  Let $F_{\infty}$ be a strongly $S$-admissible, pro-$p,$
$p$-adic  Lie extension of $F$ of dimension $>1$, and suppose that
Conjecture A holds for $Y_S(A/F^{\cyc})$. Is
$Y_S(A/F_{\infty})$ a pseudo-null $\Op\ps{G}$-module, or equivalently
a torsion $\Op\ps{H}$-module?

\medskip
When $A$ is the  (discrete) quotient module of a Galois representation coming
from arithmetic, the assertion of the question is conjectured to
always hold (see \cite{CS, Jh, LimFine, Oc}). When $A=\Qp/\Zp$, the
dual fine Selmer group is precisely
$\Gal(K(F_{\infty})/F_{\infty})$, where $K(F_{\infty})$ is the
maximal unramified pro-$p$ extension of $F_{\infty}$ at which every
prime of $F_{\infty}$ above $p$ splits completely. In this case, it
is expected that the pseudo-nullity property should hold for
admissible $p$-adic Lie extensions ``coming from algebraic
geometry'' (see \cite[Question 1.3]{HS} for details, and see also
\cite[Conjecture 7.6]{Sh1} for a related assertion and \cite{Sh2}
for positive results in this direction). However, we should mention
that Hachimori and Sharifi \cite{HS} constructed a class of
admissible $p$-adic Lie extensions $F_{\infty}$ of $F$ of dimension
$>1$ such that $\Gal(K(F_{\infty})/F_{\infty})$ is not pseudo-null.
Note that their extensions come from CM fields, where it is
generally expected that these fields cannot be realized as
admissible $p$-adic Lie extensions which are carved out by algebraic
geometrical objects (see \cite[Page 570]{HS}).

We now record a lemma which relates the module theoretical structure
of the dual fine Selmer groups and the second Iwasawa cohomology
groups over certain strongly $S$-admissible, pro-$p$, $p$-adic Lie
extensions.

\bl \label{pesudo-null H2} Let $F_{\infty}$ be a strongly
$S$-admissible, pro-$p$, $p$-adic Lie extension of $F$ of dimension
$>1$. Suppose that the following conditions are satisfied.
\begin{enumerate}
  \item[$(i)$] Conjecture A holds for $Y_S(A/F^{\cyc})$.
  \item[$(ii)$] For every $v\in S$, the decomposition group of $G$ at $v$ has dimension $\geq 2$.
       \end{enumerate}
Then $Y_S(A/F_{\infty})$ is a pseudo-null $\Op\ps{G}$-module if and
only if $\Hi^2(F_{\infty}/F, T^*)$ is a pseudo-null
$\Op\ps{G}$-module.  Furthermore, $Y_S(A/F_{\infty})(\pi)$ and
$\Hi^2(F_{\infty}/F, T^*)(\pi)$ have the same elementary
representations as $\Op\ps{H}$-modules. \el

\bpf
 By the Poitou-Tate sequence (see \cite[Eqn. (45)]{CS}), we have the following exact sequence
 \[ 0\lra Y_S(A/F_{\infty}) \lra \Hi^2(F_{\infty}/F, T^*)\lra
 \Bigg(\bigoplus_{v\in S} K_v^0(A/F_{\infty})\Bigg)^{\vee}.\]
 In view of assumption (i) and Lemma \ref{fg La H}, the first two
terms in the exact sequence are finitely generated
$\Op\ps{H}$-modules. Therefore, by virtue of Lemma
\ref{pseudo-isomorphism lemma}, the assertions of the lemma will
follow once we can show that the last term in the exact sequence is
a finitely generated torsion $\Op\ps{H}$-module with trivial
$\mu_{\Op\ps{H}}$-invariant. For each $v\in S$, fix a prime $w$ of
$F_{\infty}$ above $v$. By abuse of notation, we shall also denote
the prime of $F^{\cyc}$
  below $w$ by $w$. Denote by $H_w$ the decomposition group of $H$ at $w$. Then one
  can check that $K_v^0(A/F_{\infty})^{\vee}$ is isomorphic to a finite sum
  of terms of the form
\[  \Ind^{H}_{H_w}\Big(A(F_{\infty,w})^{\vee}\Big),
\] where $A(F_{\infty,w}) = A^{\Gal(\bar{F}_v/F_{\infty,w})}$. Since
$(A(F_{\infty,w})^{\vee}$ is finitely generated over $\Op$ and $H_w$
has dimension at least 1 by hypothesis (ii), we may apply Lemma
\ref{ind} to conclude that
$\Ind^{H}_{H_w}\left(A(F_{\infty,w})^{\vee}\right)$ is finitely
generated over $\Op\ps{H}$ with trivial $\mu_{\Op\ps{H}}$-invariant.
Thus, the lemma is proven.
 \epf

We now come to the main theme of the paper which is to study
the preservation of the pseudo-nullity property under congruences. We first
mention that under the assumption of the validity of Conjecture A,
it follows from Lemma \ref{fg La H} that $Y_S(A/F_{\infty})$ is a
finitely generated $\Op\ps{H}$-module. By a standard result of
Howson \cite[Lemma 2.7]{Ho}, this in turn implies that
$\mu_{\Op\ps{G}}(Y_S(A/F_{\infty})) =0$. By \cite[Remark 3.33]{V02},
it then follows that $Y_S(A/F_{\infty})(\pi)$ is pseudo-null as an
$\Op\ps{G}$-\textit{module}. However, being a finitely generated
$\Op\ps{H}$-module, the structure of $Y_S(A/F_{\infty})(\pi)$ is a
more subtle issue, and this is precisely the point of our next two
theorems when considering the preservation of the pseudo-nullity of
the dual fine Selmer groups under congruences.

\bt \label{ps cong high}  Let $F_{\infty}$ be a strongly
$S$-admissible pro-$p$ $p$-adic Lie extension of $F$ of dimension
$>1$. Suppose that the following conditions are satisfied.
\begin{enumerate}
 \item[$(a)$] There is an isomorphism $A[\pi^{\theta_{A}+1}]\cong B[\pi^{\theta_{A}+1}]$ of
$G_S(F)$-modules, where $\theta_A :=
\theta_{\Op\ps{H}}\big(Y_S(A/F_{\infty})\big)$.
 \item[$(b)$] Conjecture A holds for $Y_S(A/F^{\cyc})$ $($and hence for $Y_S(B/F^{\cyc})$ by Proposition \ref{fg cong} $)$.
   \item[$(c)$] $Y_S(A/F_{\infty})$ is a pseudo-null $\Op\ps{G}$-module.
  \item[$(d)$] For each $v\in S$, the decomposition group of $G$ at $v$ has dimension $\geq 2$.
    \end{enumerate}
    Then $Y_S(B/F_{\infty})$ is a pseudo-null $\Op\ps{G}$-module. Furthermore, $Y_S(A/F_{\infty})(\pi)$ and $Y_S(B/F_{\infty})(\pi)$ have the same elementary representations as $\Op\ps{H}$-modules. \et

\br
 Note that it follows from assumption (b)
that $Y_S(A/F_{\infty})(\pi)$ and $Y_S(B/F_{\infty})(\pi)$ have the
same elementary representations as $\Op\ps{G}$-modules. In fact, as
seen in the dicussion before Theorem \ref{ps cong high}, they are
pseudo-null as $\Op\ps{G}$-modules. \textbf{The point of our theorem
is that we can relate the $\Op\ps{H}$-module structure of the
$\pi$-primary submodules of the two dual fine Selmer groups}. \er

\bpf[Proof of Theorem \ref{ps cong high}] For every $i\geq 0$ and
every open subgroup $G_0$ of $G$, one has an isomorphism
\[ \Ext_{\Op\ps{G}}^{i}(M, \Op\ps{G}) \cong
\Ext_{\Op\ps{G_0}}^{i}(M, \Op\ps{G_0})\] for any finitely generated
$\Op\ps{G}$-module $M$ (cf.\ \cite[Proposition 5.4.17]{NSW}).
Therefore, replacing $F$ if necessary, we may assume that $H$ and
$G$ are uniform pro-$p$ groups. By Lemma \ref{pesudo-null H2}, we
may work with second Iwasawa cohomology groups which we shall do for
the proof of the theorem. By Proposition \ref{pseudo-isomorphic2},
it remains to show that
 \[ \mu_{\Op\ps{H}}\Big(\Hi^2(F_{\infty}/F, T_A^*)/\pi^{n}\Big) =
\mu_{\Op\ps{H}}\Big(\Hi^2(F_{\infty}/F, T_B^*)/\pi^{n}\Big)\] for
$1\leq n \leq \theta_H(A)+1$. Fix such an arbitrary $n$. We now note
that it follows from \cite[Lemma 2.1]{LimFine} that
$\Hi^2(F_{\infty}/F, T_Z^*)/\pi^n\cong \Hi^2(F_{\infty}/F,
T_Z^*/\pi^n T_Z^*)$ for $Z=A,B$.
% Changes made on 4/2/2017
By the congruence condition (a) of the hypotheses, we
have an isomorphism $T_A^*/\pi^n T_A^*\cong T_B^*/\pi^n T_B^*$ of
$G_S(F)$-modules which in turn induces an isomorphism
\[ \Hi^2(F_{\infty}/F, T_A^*/\pi^n
T_A^*) \cong \Hi^2(F_{\infty}/F, T_B^*/\pi^n T_B^*) \] of
$\Op\ps{H}$-modules. The required equality of
$\mu_{\Op\ps{H}}$-invariants is now an immediate consequence of
this. Thus, we have proven the theorem. \epf

When $\mu_{\Op\ps{H}}(Y_S(A/F_{\infty})) =0$, we can prove the
following theorem which has slightly weaker requirement on the
decomposition groups of the $p$-adic Lie extension. Due to this
weaker assumption, we cannot employ Lemma \ref{pesudo-null H2}, and
so we have to work with the dual fine Selmer groups directly.

\bt \label{ps cong}  Let $F_{\infty}$ be a strongly $S$-admissible
pro-$p$ $p$-adic Lie extension of $F$ of dimension $>1$. Suppose
that the following conditions are satisfied.
\begin{enumerate}
 \item[$(a)$] There is an isomorphism $A[\pi]\cong B[\pi]$ of
$G_S(F)$-modules.
 \item[$(b)$] Conjecture A holds for $Y_S(A/F^{\cyc})$ $($and hence for $Y_S(B/F^{\cyc})$ by Proposition \ref{fg cong} $)$.
  \item[$(c)$] We have $\mu_{\Op\ps{H}}(Y_S(A/F_{\infty})) =0$,
  and hence $Y_S(A/F_{\infty})$ is a pseudo-null $\Op\ps{G}$-module.
    \item[$(d)$] For each $v\in S$, either one of the following holds.
   \begin{enumerate}
    \item[$(i)$] The decomposition group of $G$ at $v$, denoted by $G_v$, has dimension $\geq 2$.
    \item[$(ii)$] For every prime $w$ of $F_{\infty}$ above $v$, $B(F_{\infty ,w})$ is a divisible $\Op$-module.
 \end{enumerate}
    \end{enumerate}
    Then $Y_S(B/F_{\infty})$ is a pseudo-null $\Op\ps{G}$-module. Moreover, we have
    $\mu_{\Op\ps{H}}(Y_S(B/F_{\infty})) =0$. \et

\bpf
  It follows from Lemma \ref{fg La H} and assumptions (b) and (c) that
$Y_S(A/F_{\infty})$ is a finitely generated torsion
$\Op\ps{H}$-module. By \cite[Corollary 1.10]{Ho} (or
\cite[Proposition 4.12]{LimFine}), we then have \[ 0 =
\rank_{\Op\ps{H}}\big(Y_S(A/F_{\infty})\big) =
\rank_{k\ps{H}}\big(Y_S(A/F_{\infty})/\pi\big) -
\rank_{k\ps{H}}\big(Y_S(A/F_{\infty})[\pi]\big).\] On the other
hand, by an application of \cite[Proposition 1.6, Corollary
1.7]{Ho}, we have
\[ \rank_{k\ps{H}}\big(Y_S(A/F_{\infty})[\pi]\big)= \mu_{\Op\ps{H}}\big(Y_S(A/F_{\infty})[\pi]\big),\]
and the latter quantity is zero as a consequence of assumption (c).
Hence we may conclude that $Y_S(A/F_{\infty})/\pi$ is a torsion
$k\ps{H}$-module. To continue, we need to introduce the $\pi$-fine
Selmer groups. For a prime $v$ in $S$, and for each finite extension
$L$ of $F$ contained in $F_{\infty}$, we set
 \[K_v^1(A[\pi]/L) = \bigoplus_{w|v}H^1(L_w, A[\pi]),\]
where $w$ runs over the (finite) set of primes of $L$ above $v$. The
$\pi$-fine Selmer group (with respect to $S$) is then defined to be
\[ R_S(A[\pi]/L) = \ker\Big(H^1(G_S(L),A[\pi])
\longrightarrow \bigoplus_{v\in S}K_v^1(A[\pi]/L)\Big). \] Set
$R_S(A[\pi]/F_{\infty})=\ilim_L R_S(A[\pi]/L)$, where $L$ runs
through all finite subextensions of $F_{\infty}/F$. We then have the
following diagram
\[  \entrymodifiers={!! <0pt, .8ex>+} \SelectTips{eu}{}\xymatrix{
    0 \ar[r]^{} & R_S(A[\pi]/F_{\infty}) \ar[d]_{f_A} \ar[r] &
    H^1(G_S(F_{\infty}), A[\pi])
    \ar[d]_{g_A}
    \ar[r]^{} & \bigoplus_{v\in S}K_v(A[\pi]/F_{\infty}) \ar[d]_{h_A} \\
    0 \ar[r]^{} & R_S(A/F_{\infty})[\pi] \ar[r]^{}
    & H^1(G_S(F_{\infty}), A)[\pi] \ar[r] & \
    \bigoplus_{_{v\in S}}K_v(A/F_{\infty})[\pi]  } \]
with exact rows. The long exact sequence in cohomology arising from
$0\lra A[\pi]\lra A\lra A\lra 0$ shows that $\ker g_A$ is finite.
Thus, $\ker f_A$ is also finite. Since $H$ has dimension $\geq 1$,
the Pontryagin dual of $\ker f_A$ is a torsion $k\ps{H}$-module.
Therefore, we have a $k\ps{H}$-homomorphism
\[ Y_S(A/F_{\infty})/\pi\lra Y_S(A[\pi]/F_{\infty}) \]
with cokernel which is a torsion $k\ps{H}$-module. Since we have
already shown above that $Y_S(A/F_{\infty})/\pi$ is a torsion
$k\ps{H}$-module, it follows that $Y_S(A[\pi]/F_{\infty})$ is also a
torsion $k\ps{H}$-module. By assumption (a), this in turn implies
that $Y_S(B[\pi]/F_{\infty})$ is a torsion $k\ps{H}$-module. Now
consider the following diagram
\[  \entrymodifiers={!! <0pt, .8ex>+} \SelectTips{eu}{}\xymatrix{
    0 \ar[r]^{} & R_S(B[\pi]/F_{\infty}) \ar[d]_{f_B} \ar[r] &
    H^1(G_S(F_{\infty}), B[\pi])
    \ar[d]_{g_B}
    \ar[r]^{} & \bigoplus_{v\in S}K_v(B[\pi]/F_{\infty}) \ar[d]_{h_B} \\
    0 \ar[r]^{} & R_S(B/F_{\infty})[\pi] \ar[r]^{}
    & H^1(G_S(F_{\infty}), B)[\pi] \ar[r] & \
    \bigoplus_{_{v\in S}}K_v(B/F_{\infty})[\pi]  } \]
with exact rows. By a similar cohomological argument as above, the
map $g_B$ has finite kernel and trivial cokernel. We shall now show
that $\ker h_B$ is a cofinitely generated torsion $k\ps{H}$-module.
Write $h_B =\oplus_w h_{B, w}$, where $w$ runs over the (finite) set
of primes of $F^{\cyc}$ above $S$. Denote by $H_w$ the decomposition
group of $F_{\infty}/F^{\cyc}$ corresponding to a fixed prime of
$F_{\infty}$, which we also denote by $w$, above $w$. Then we have
$\ker h_{B,w}=\mathrm{Coind}^{H_w}_H\big(B(F_{\infty,w})/\pi\big)$,
where $v$ is the prime of $F$ below $w$. If $v$ satisfies assumption
(d)(ii), then $\ker h_{B,w}= 0$. Now suppose that $v$ satisfies
assumption (d)(i), a similar argument to that in Lemma
\ref{pesudo-null H2} shows that $\ker h_{B,w}$ is a cotorsion
$k\ps{H}$-module. In conclusion, we have shown that $\ker h_B$ is a
cofinitely generated torsion $k\ps{H}$-module. By a diagram chasing
argument, one then has that the cokernel of the map $f_B$ is a cofinitely
generated torsion $k\ps{H}$-module. Hence we have a
 $k\ps{H}$-homomorphism
\[ Y_S(B/F_{\infty})/\pi\lra Y_S(B[\pi]/F_{\infty}) \]
whose kernel and cokernel are torsion $k\ps{H}$-modules. Combining
this with the above observation that $Y_S(B[\pi]/F_{\infty})$ is
torsion over $k\ps{H}$, we have that $Y_S(B/F_{\infty})/\pi$ is
torsion over $k\ps{H}$. By \cite[Corollary 1.10]{Ho} (or
\cite[Proposition 4.12]{LimFine}), we have
\[ \rank_{\Op\ps{H}}\big(Y_S(B/F_{\infty})\big) =
\rank_{k\ps{H}}\big(Y_S(B/F_{\infty})/\pi\big) -
\rank_{k\ps{H}}\big(Y_S(B/F_{\infty})[\pi]\big).\] From these, it
follows that $Y_S(B/F_{\infty})$ is a finitely generated torsion
$\Op\ps{H}$-module and $Y_S(B/F_{\infty})[\pi]$ is a finitely
generated torsion $k\ps{H}$-module. The latter is equivalent to
$\mu_{\Op\ps{H}}\big(Y_S(B/F_{\infty})\big) =0$ by \cite[Remark
3.33]{V02}. \epf

\section{Comparing Galois groups and fine Selmer groups} \label{Galois and fine}

As before, $p$ denote a prime and $F$ a number field. If $p=2$,
assume further that $F$ has no real primes. Let $\Op$ be the ring of
integers of a fixed finite extension $K$ of $\Qp$. Let $A$ denote
the quotient module of a finite dimensional $K$-vector space, which
is endowed with a continuous $G_S(F)$-action for a finite set $S$ of
primes. We shall also assume that the set $S$ contains all the
primes above $p$, the ramified primes of $A$ and the infinite
primes. Inspired by the relation between the Iwasawa
$\mu$-conjecture and Conjecture A, the first author is led to ask
the following question (see \cite{LimPS}).

\medskip

\noindent \textbf{Question B$'$:}  Let $F_{\infty}$ be an
$S$-admissible $p$-adic Lie extension of a number field $F$ of
dimension $>1$ with the property that $G_S(F_{\infty})$ acts
trivially on $A[\pi]$. Suppose that $\Gal(K(F_{\infty})/F_{\infty})$
is a finitely generated $\Zp\ps{H}$-module. Can one deduce that
$Y_S(A/F_{\infty})$ is a pseudo-null $\Op\ps{G}$-module from the
knowledge that $\Gal(K(F_{\infty})/F_{\infty})$ is a pseudo-null
$\Zp\ps{G}$-module?

\medskip
In the same paper, the first author gave a partial answer to the
above question (see \cite[Theorem 2.3]{LimPS}). We now apply the
criterion in Section \ref{fine Selmer group} to derive refinements
of this result. As a start, we have the following.

 \bp \label{pseudo-null main}
 Suppose that $F$ contains $\mu_p$. Let $F_{\infty}$ be a strongly $S$-admissible pro-$p$
$p$-adic Lie extension of $F$ of dimension $>1$.
Suppose that the following conditions are satisfied.
\begin{enumerate}
 \item[$(a)$] $G_S(F_{\infty})$ acts trivially on $A[\pi]$.
 \item[$(b)$] Conjecture A holds for $\Gal(K(F_{\infty})/F_{\infty})$.
  \item[$(c)$] $\mu_{\Zp\ps{H}}(\Gal(K(F_{\infty})/F_{\infty})) =0$. In particular, $\Gal(K(F_{\infty})/F_{\infty})$ is a pseudo-null $\Zp\ps{G}$-module.
  \item[$(d)$] For each $v\in S$, either one of the following holds.
   \begin{enumerate}
    \item[$(i)$] The decomposition group of $G$ at $v$, denoted by $G_v$, has dimension $\geq 2$.
    \item[$(ii)$] For every prime $w$ of $F_{\infty}$ above $v$, $A(F_{\infty ,w})$ is a divisible $\Op$-module.
 \end{enumerate}
    \end{enumerate}
    Then $Y_S(A/F_{\infty})$ is a
pseudo-null $\Op\ps{G}$-module. \ep

\bpf
Since $\Op$ is free over $\Zp$, we have an isomorphism
  \[ \Ext_{\Op\ps{G}}^{i}(M, \Op\ps{G}) \cong \Op\ot_{\Zp}
\Ext_{\Zp\ps{G}}^{i}(M, \Zp\ps{G})\]
  of $\Op\ps{G}$-modules for every $\Op\ps{G}$-module $M$. Therefore, it suffices to show that
$Y_S(A/F_{\infty})$ is a pseudo-null $\Zp\ps{G}$-module.
  Replacing $F$ by a larger extension if necessary, we may assume that
$A[\pi]$ is a trivial $G_S(F)$-module. It then follows that
$A[\pi]\cong (\Op/\pi)^d \cong (\Z/p)^{fd}$ as $G_S(F)$-modules.
Here $d=\mathrm{corank_{\Op}}(A)$ and $f=[K:\Qp]$. Since
$Y_S\big((\Qp/\Zp)^{fd}/F_{\infty}\big) =
\Gal(K(F_{\infty}/F_{\infty})^{fd}$, it follows from assumption (c)
that
$\mu_{\Zp\ps{H}}\Big(Y_S\big((\Qp/\Zp)^{fd}/F_{\infty}\big)\Big)=0$
and $Y_S\big((\Qp/\Zp)^{fd}/F_{\infty}\big)$ is a pseudo-null
$\Zp\ps{G}$-module. The conclusion of the proposition now follows
from an application of Theorem \ref{ps cong}. \epf

The next result considers the case when
$\mu_{\Zp\ps{H}}(\Gal(K(F_{\infty})/F_{\infty})) \neq0$. This result
also gives a relation between the $\pi$-primary submodule of the
dual fine Selmer group and the $p$-primary submodule of
$\Gal(K(F_{\infty})/F_{\infty})$.

\bp \label{ps cong high Galois}  Let $F_{\infty}$ be a strongly
$S$-admissible pro-$p$ $p$-adic Lie extension of $F$ of dimension
$>1$.  Suppose that the following conditions are satisfied.
\begin{enumerate}
 \item[$(a)$] $\Gal(K(F_{\infty})/F_{\infty})$ is a finitely generated
torsion $\Zp\ps{H}$-module.
 \item[$(b)$] $G_S(F_{\infty})$ acts trivially on $A[\pi^{n}]$, where
\[ n = \Bigg\lceil \frac{\theta_{\Op\ps{H}}\big(\Gal(K(F_{\infty})/F_{\infty})\big)+1}{e} \Bigg\rceil \]
and $e$ is the ramification index of $\Op/\Zp$.
    \item[$(c)$] For each $v\in S$, the decomposition group of $G$ at $v$ has dimension $\geq 2$.
    \end{enumerate}
Then $Y_S(A/F_{\infty})$ is a finitely generated torsion
$\Op\ps{H}$-module, and hence a pseudo-null $\Op\ps{G}$-module.
Furthermore, $Y_S(A/F_{\infty})(\pi)$ and
$\Gal(K(F_{\infty})/F_{\infty})(p)^{fd}$ have the same elementary
representations as $\Zp\ps{H}$-modules, and we have
\[ e\theta_{\Op\ps{H}}\big(Y_S(A/F_{\infty})\big) =  \theta_{\Zp\ps{H}}\big(\Gal(K(F_{\infty})/F_{\infty})\big).\]  Here $d=\mathrm{corank_{\Op}}(A)$ and $f=[K:\Qp]$.
\ep

\bpf[Proof of Proposition \ref{ps cong high Galois}]
 Via a similar observation to that in Proposition \ref{pseudo-null main},
it suffices to show that $Y_S(A/F_{\infty})$ is a pseudo-null
$\Zp\ps{G}$-module. Also, replacing $F$ if necessary, we may assume
that $A[\pi^{n}] = A[p^{en}]$ is a trivial $G_S(F)$-module. It then
follows that $A[p^{en}]\cong B[p^{en}]$ as $G_S(F)$-modules, where
$B = (\Qp/\Zp)^{fd}$. By assumption (b), we have that $ne \geq
\theta_{\Op\ps{H}}\big(\Gal(K(F_{\infty})/F_{\infty})\big)+1$. In
view of assumption (a), we may apply Theorem \ref{ps cong high} to
obtain the conclusion of the proposition. \epf

% Changes made and added on 20/3/2017

\section{Some further remarks} \label{further remark section}

In this section, we mention a result which gives a relation between Conjecture A and the structure
of the Selmer group of the residual representation. In the process of obtaining such a relation, we need to compare the Selmer group and the so-called strict Selmer group of Greenberg \cite{G89}.
To prepare for this, we recall certain background material from \cite{CG}. Let $\mathcal{K}$ be a finite extension of $\Qp$ and $E$ an elliptic curve defined over $\mathcal{K}$. Throughout our discussion, we shall always assume that our elliptic curve $E$ has good ordinary reduction. For every algebraic extension $\mathcal{L}$ of $\mathcal{K}$, we denote by $\m_\mathcal{L}$ the maximal ideal of the ring of integers of $\mathcal{L}$. Let $\hat{E}(\m_{\mathcal{L}})$ be the formal group of $E$ over the ring of integers of $\mathcal{L}$. Kummer theory gives us the following commutative diagram
\[  \entrymodifiers={!! <0pt, .8ex>+} \SelectTips{eu}{}\xymatrix{
    0 \ar[r]^{} & \hat{E}(\m_{\mathcal{K}})/p \ar[d] \ar[r] &
    H^1(\mathcal{K}, \hat{E}(\m_{\overline{\mathcal{K}}})[p])
    \ar[d]
    \ar[r]^{} & H^1(\mathcal{K}, \hat{E}(\m_{\overline{\mathcal{K}}})[p]\ar[r] \ar[d]&0 \\
    0 \ar[r] & E(\mathcal{K})/p \ar[r]^{}
    & H^1(\mathcal{K}, E[p]) \ar[r] & \
    H^1(\mathcal{K}, E)[p] \ar[r] &0 } \]
with exact rows. Let $\mathcal{K}_{\infty}$ be the cyclotomic $\Zp$-extension of $\mathcal{K}$. Then we have a similar commutative diagram as above for each intermediate field of the  extension $\mathcal{K}_{\infty}/\mathcal{K}$. Taking direct limit of these diagrams and noting that $H^1(\mathcal{K}_{\infty},\hat{E}(\m_{\overline{\mathcal{K}}}))=0$ (cf.\ \cite[Corollary 3.2]{CG}), we obtain the following commutative diagram
\[  \entrymodifiers={!! <0pt, .8ex>+} \SelectTips{eu}{}\xymatrix{
     & \hat{E}(\m_{\mathcal{K}_{\infty}})/p \ar[d] \ar[r]^(.4){\cong} &
    H^1(\mathcal{K}_{\infty}, \hat{E}(\m_{\overline{\mathcal{K}}})[p])
    \ar[d]^{\lambda_p}
    & & \\
    0 \ar[r] & E(\mathcal{K}_{\infty})/p \ar[r]^{\kappa_p}
    & H^1(\mathcal{K}_{\infty}, E[p]) \ar[r] & \
    H^1(\mathcal{K_{\infty}}, E)[p] \ar[r] &0 } \]
with exact row.

\bl \label{local compare}
With notation as above, we have a surjection
\[ H^1(\mathcal{K}_{\infty},\tilde{E}[p])\twoheadrightarrow H^1(\mathcal{K}_{\infty},E)[p], \]
where $\tilde{E}$ denotes the reduction of $E$ mod $v$.
\el

\bpf
 From the commutative diagram before the lemma, we have that $\im\la_p\subseteq \im \kappa_p$ and this in turn induces a surjection
\[ H^1(\mathcal{K}_{\infty}, E[p])/\im\la_p \twoheadrightarrow H^1(\mathcal{K}_{\infty}, E[p])/\im\kappa\cong  H^1(\mathcal{K_{\infty}}, E)[p]. \]
It remains to show that $ H^1(\mathcal{K}_{\infty}, E[p])/\im\la_p\cong H^1(\mathcal{K}_{\infty},\tilde{E}[p])$.
But this is a consequence of the long exact sequence in cohomology arising from
$$
0\lra \hat{E}(\m_{\overline{\mathcal{K}}})[p]\lra E[p]\lra \tilde{E}[p]\lra 0
$$
which also  gives the following exact sequence
\[H^1(\mathcal{K}_{\infty}, \hat{E}(\m_{\overline{\mathcal{K}}})[p])\stackrel{\la_p}{\lra}
 H^1(\mathcal{K}_{\infty}, E[p]) \lra H^1(\mathcal{K}_{\infty}, \tilde{E}[p]) \lra H^2(\mathcal{K}_{\infty}, \hat{E}(\m_{\overline{\mathcal{K}}})[p]),\]
 and noting that the last term is zero by the fact that $\Gal(\overline{\mathcal{K}}/\mathcal{K}_{\infty})$ has $p$-cohomological dimension $\leq 1$ (cf. \cite[Theorem 7.1.8(i)]{NSW}).
\epf

We now turn to the number field context. Let $F$ be a number field. From now on, $E$ will denote an elliptic curve defined over $F$.
Throughout our discussion, we always assume that $E$ has good ordinary reduction
at all primes of $F$ above $p$. The classical $p^n$-Selmer group of $E$ is defined to be
\[ S(E[p^n]/F) = \ker\Big(H^1(F, E[{p^n}])\lra \bigoplus_{v} H^1(F_v, E)[p^n]
\Big), \]
where $v$ runs through all the primes of $F$. Let $S$ be
a finite set of primes of $F$ which contains the primes above $p$,
the bad reduction primes of $E$ and the infinite primes. By \cite[p. 8]{CS10}, we have the following equivalent description of the $p^n$-Selmer group
\[ S(E[p^n]/F) = \ker\Big(H^1(G_S(F), E[{p^n}])\lra \bigoplus_{v\in S} H^1(F_v, E)[p^n]
\Big).\]
We set $S(E/F)= \ilim_n S(E[p^n]/F)$ which is precisely the classical $p$-primary Selmer group.
Let $F_{\infty}$ be a strongly admissible pro-$p$, $p$-adic Lie extension of $F$. Write $G= \Gal(F_{\infty}/F)$ and $H= \Gal(F_{\infty}/F^{\cyc})$. Similarly, we can define $S(E[p^n]/L)$ and $S(E/L)$ for each finite extension $L$ of $F$ which is contained in $F_{\infty}$. We then set $S(E[p^n]/F_{\infty})=\ilim_L S(E[p^n]/L)$ and $S(E/F_{\infty})=\ilim_L S(E/L)$. Writing $S(F_{\infty})$ for the set of primes of $F_{\infty}$ above $S$, it is not difficult to verify that
\[ S(E[p^n]/F_{\infty}) = \ker\Big(H^1(G_S(F_{\infty}), E[{p^n}])\lra \bigoplus_{w\in S(F_{\infty})} H^1(F_{\infty, w}, E)[p^n]
\Big) \]
and
\[ S(E/F_{\infty}) = \ker\Big(H^1(G_S(F_{\infty}), E[{p^{\infty}}])\lra \bigoplus_{w\in S(F_{\infty})} H^1(F_{\infty,w}, E)[p^{\infty}]
\Big). \]
 Furthermore, we have that
$S(E[p^n]/F_{\infty})= \ilim_L S(E[p^n]/L^{\cyc})$ and $S(E/F_{\infty})=\ilim_L S(E/L^{\cyc})$.

We now introduce the strict Selmer group of Greenberg \cite{G89}. For each $v\in S$, we define the strict $p^n$-Selmer group of $E$ over $F$ to be
\[ S^{str}(E[p^n]/F) = \ker\Big(H^1(F, E[{p^n}])\lra \bigoplus_{v\in S} H^1(F_v,D_v[p^n])
\Big), \]
 where $D_v$ is taken to be $E[p^{\infty}]$ or $\widetilde{E}_v[p^{\infty}]$
according as $v$ does not or does divide $p$, and here $\widetilde{E}_v$ denotes the
reduction of $E$ mod $v$. We then set $S^{str}(E[p^n]/F)=\ilim S^{str}(E[p^n]/F$). For each intermediate subfield $L$ of $F_{\infty}/F$, we have analogous definition for  $S^{str}(E[p^n]/L)$, and we set $S^{str}(E[p^n]/F_{\infty}) = \ilim_L S^{str}(E[p^n]/L)$. The strict ($p^{\infty}$-)Selmer group of $E$ over $F_{\infty}$ is then given by $\ilim_n S^{str}(E[p^n]/F_{\infty})$.  It is a straightforward exercise to verify that
\[ S^{str}(E[p^n]/F_{\infty}) = \ker\Big(H^1(G_S(F_{\infty}), E[{p^n}])\lra \bigoplus_{w\in S(F_{\infty})} H^1(F_{\infty,w}, D_w[p^n])
\Big) \]
and
\[ S^{str}(E/F_{\infty}) = \ker\Big(H^1(G_S(F_{\infty}), E[{p^{\infty}}])\lra \bigoplus_{w\in S(F_{\infty})} H^1(F_{\infty,w}, D_w[p^{\infty})
\Big). \]
As before, we also have
$S^{str}(E[p^n]/F_{\infty})= \ilim_L S^{str}(E[p^n]/L^{\cyc})$ and $S^{str}(E/F_{\infty})=\ilim_L S^{str}(E/L^{\cyc})$.

It is well-known that $H^1(F_{\infty,w}, D_w[p^{\infty}]) = H^1(F_{\infty,w}, E)[p^{\infty}]$. (This can be easily verified when $w$ does not divide p; in the event that $w$ divides $p$, this follows from \cite[Proposition 4.8]{CG}.) As a consequence, we have $S(E/F_{\infty}) = S^{str}(E/F_{\infty})$. Write $X(E/F_{\infty})$ for the Pontryagin dual of $S(E/F_{\infty})$. On the other hand, for the $p$-Selmer groups, we have the following commutative diagram
\[  \entrymodifiers={!! <0pt, .8ex>+} \SelectTips{eu}{}\xymatrix{
    0 \ar[r]^{} & \Sel^{str}(E[p]/F_{\infty}) \ar[d]_{a} \ar[r] &
    H^1(G_S(F_{\infty}), E[p])
    \ar@{=}[d]
    \ar[r]^(.45){\psi_s} & \displaystyle\bigoplus_{w\in S(F_{\infty})} H^1(F_{\infty,w}, D_w[p]) \ar[d]_(.6){\varphi} \\
    0 \ar[r]^{} & \Sel(E[p]/F_{\infty}) \ar[r]^{}
    & H^1(G_S(F_{\infty}), E[p]) \ar[r]^(.45){\psi} & \
    \displaystyle \bigoplus_{w\in S(F_{\infty})} H^1(F_{\infty,w}, E)[p]  } \]
     with exact rows and the map $\varphi$ is surjective.

We can now state the main result of this section.

\bt \label{Sel vs fine} Let $E$ be an elliptic curve defined over a number field $F$
which has good ordinary reduction at all primes above $p$. Let $F_{\infty}$ be a strongly admissible pro-$p$, $p$-adic Lie extension of $F$.
Suppose that $X(E/F_{\infty})$ is torsion over $\Zp\ps{G}$.
Then the following statements are equivalent.

\begin{enumerate}
         \item[$(1)$] $X(E/F_{\infty})$ is finitely generated over $\Zp\ps{H}$.
  \item[$(2)$] We have $H^2(G_S(F_{\infty}), E[p])=0$ and there is a short exact sequence
\[0 \lra S^{str}(E[p]/F_{\infty}) \lra H^1(G_S(F_{\infty}), E[p]) \lra \bigoplus_{w\in S(F_{\infty})} H^1(F_{\infty,w}, D_w[p])\lra 0.\]
\end{enumerate}
\et

\br \label{Su remark} Note that the assertion
``$H^2(G_S(F_{\infty}), E[p])=0$" in statement (2) is equivalent to
Conjecture A being valid for $Y(E/F^{\cyc})$ (see \cite[Proposition
4.6]{Su}). (Actually, the said assertion is only proved in
\cite[Proposition 4.6]{Su} for the cyclotomic situation but it is
not difficult to check that the same argument carries over to the
general case.)  \er

Theorem \ref{Sel vs fine} can be viewed as the mod $p$ analogue of the following result, namely: the dual Selmer over an admissible $p$-adic Lie extension is torsion if and only if the defining sequence for the Selmer group is short exact and $H^2(G_S(F_{\infty}),E[p^{\infty}])=0$..  We now apply Theorem \ref{Sel vs fine} to give another different, and slightly more conceptual, proof of a result of
Vatsal and Greenberg (see \cite[Page 18, Statement A]{GV}).

\bc \label{GV proof} Let $E$ and $E'$ be two elliptic curves defined over a number field $F$
which have good ordinary reduction at all primes above $p$. Suppose that $E[p]\cong E'[p]$ as $\Gal(\bar{F}/F)$-modules. Let $F_{\infty}$ be a strongly admissible pro-$p$, $p$-adic Lie extension of $F$. Then $X(E/F_{\infty})$ is finitely generated over $\Zp\ps{H}$ if and only if $X(E'/F_{\infty})$ is finitely generated over $\Zp\ps{H}$.
\ec

\bpf
 The corollary follows from combining Theorem \ref{Sel vs fine} with the hypothesis that $E[p]\cong E'[p]$. \epf

 We record another interesting corollary of Theorem \ref{Sel vs fine}. Namely, we give a sufficient condition which ensures that  the defining sequence for the classical Selmer group of $E[p]$ is surjective. 

\bc Let $E$ be an elliptic curve defined over a number field $F$
which has good ordinary reduction at all primes above $p$.  Let $F_{\infty}$ be a strongly admissible pro-$p$, $p$-adic Lie extension of $F$.
Suppose that $X(E/F^{\cyc})$ is finitely generated over $\Zp$.
Then we have a short exact sequence
\[0 \lra S(E[p]/F_{\infty}) \lra H^1(G_S(F_{\infty}), E[p]) \lra \bigoplus_{w\in S(F_{\infty})} H^1(F_{\infty,w}, E)[p]\lra 0.\]
\ec

\bpf
 This follows from Theorem \ref{Sel vs fine} and the diagram before Theorem \ref{Sel vs fine}.
 By \textit{loc. cit.}, the map $\psi_s$ is surjective. It then follows from the commutativity of the rightmost square that $\psi$ is also surjective which gives the required short exact sequence.
\epf

It therefore remains to prove Theorem \ref{Sel vs fine}. We shall first require a lemma. Set $S^{*}(E[p]/F_{\infty}) = \plim_L \big(S^{str}(E[p]/L)\big)$.

\bl \label{inject} We have an injection
 \[S^*(E[p]/F_{\infty})\hookrightarrow\Hom_{\mathbb{F}_p\ps{G}}(S^{str}(E[p]/F_{\infty})^{\vee},\mathbb{F}_p\ps{G}). \]
 In particular, if $S^{str}(E[p]/F_{\infty})^{\vee}$ is torsion over $\mathbb{F}_p\ps{G}$, then $S^*(E[p]/F_{\infty})=0$.
\el

\bpf
   The proof is quite similar to that in \cite[Proposition 7.1]{HV}. For the convenience of the reader, we shall give a detailed proof here.
   For each finite extension $L$ of $F$ contained in $F_{\infty}$, we write $G_L=\Gal(F_{\infty}/L)$. The restriction maps on cohomology induces a map on Selmer groups
  \[ r_L:S^{str}(E[p]/L)\lra S^{str}(E[p]/F_{\infty})^{G_L}, \]
  whose kernel is contained in $H^1(G_L, E(F_{\infty})[p])$.
   From this, we have an exact sequence
  \[ 0\lra\plim_m\ker(r_L) \lra
 S^*(E[p]/F_{\infty}) \lra \plim_L S^{str}(E[p]/L)^{G_L}.\]
  Here the inverse limit is taken with respect to corestriction for the second term, and for the last term, the
inverse limit is taken with respect to the map induced by the following map
\[S^{str}(E[p]/F_{\infty})^{G_{L'}} \lra S^{str}(E[p]/F_{\infty})^{G_{L}},  \quad x\mapsto \sum_{\sigma\in
\Gal(L'/L)}\sigma(x)\] for $L'\supseteq L$. Now for sufficiently
large enough $L$, we have $H^1(G_L, E(F_{\infty})[p])=
E(F_{\infty})[p]$. Therefore, for these $L$, the group $G_L$ acts
trivially on $H^1(G_L, E(F_{\infty})[p])$ and hence $\ker (r_L)$. It
then follows that the corestriction map from $\ker (r_{L'})$ to
$\ker(r_L)$ is precisely multiplication by $[L':L]$ for $L'\supseteq
L$. Since these groups are killed by $p$, the said map is the zero
map. Thus, we have $\plim_L\ker (r_L) =0$. Consequently, we have an
injection
\[ S^*(E[p]/F^{\cyc}) \hookrightarrow \plim_L S^{str}(E[p]/F_{\infty})^{G_L}. \]

Now if $M$ is a finite group killed by $p$, we have a natural non-degenerate pairing $M\times M^{\vee}\lra \mathbb{F}_p$ which induces an isomorphism
   \[ M\stackrel{\cong}{\lra}\Hom_{\mathbb{F}_p}(M^{\vee},\mathbb{F}_p). \]
   Applying this to each $S^{str}(E[p]/F_{\infty})^{G_L}$ and taking inverse limit, we have
   \[ \plim_L S^{str}(E[p]/F_{\infty})^{G_L} \cong \plim_L \Hom_{\mathbb{F}_p}\left(\big(S^{str}(E[p]/F_{\infty})^{\vee}\big)_{G_L},\mathbb{F}_p\right).\]
  On the other hand, for a finitely generated $\mathbb{F}_p\ps{G}$-module $N$, we have an isomorphism
\[ \Hom_{\mathbb{F}_p}(N_{G_L},\mathbb{F}_p)\cong \Hom_{\mathbb{F}_p\ps{G}}(N,\mathbb{F}_p[G/G_L]). \]
Applying this isomorphism to each
$\big(S^{str}(E[p]/F_{\infty})^{\vee}\big)_{G_L}$ and taking limit,
we have
\[ \plim_L \Hom_{\mathbb{F}_p}\left(\big(S^{str}(E[p]/F_{\infty})^{\vee}\big)_{G_L},\mathbb{F}_p\right)\cong
\plim_m\Hom_{\mathbb{F}_p\ps{G}}\left(S^{str}(E[p]/F_{\infty})^{\vee},\mathbb{F}_p[G/G_L]\right)
\]
\[\hspace{1.8in} \cong \Hom_{\mathbb{F}_p\ps{G}}\big(S^{str}(E[p]/F_{\infty})^{\vee},\mathbb{F}_p\ps{G}\big).\]
Combining these observations, we obtain the required injection. \epf

\br
 When $F_{\infty}=F^{\cyc}$, one can even show that the injection in
 Lemma \ref{inject} is an isomorphism.
\er

We can now give the proof of Theorem \ref{Sel vs fine}.

\bpf[Proof of Theorem \ref{Sel vs fine}] Suppose that statement (1)
holds. Being a quotient of $X(E/F_{\infty})$, $Y(E/F_{\infty})$ is
therefore finitely generated over $\Z_p\ps{H}$. By \cite[Proposition
4.6]{Su} (and see Remark \ref{Su remark}), we then have
$H^2(G_S(F_{\infty}), E[p])=0$. In view of this, the Poitou-Tate
sequence then gives us the following exact sequence
\[0 \lra S^{str}(E[p]/F_{\infty}) \lra H^1(G_S(F_{\infty}), E[p]) \lra\bigoplus_{w\in S(F_{\infty})}
 H^1(F_{\infty,w}, D_w[p])\lra S^*(E[p]/F_{\infty})^{\vee}\lra 0.\]
By virtue of statement (1), $S(E[p]/F_{\infty})^{\vee}$ is finitely
generated over $\mathbb{F}_p\ps{H}$, and hence torsion over
$\mathbb{F}_p\ps{G}$. By Lemma \ref{inject}, this implies that
$S^*(E[p]/F_{\infty})=0$ which in turn gives the required short
exact sequence.

Conversely, suppose that $ H^2(G_S(F_{\infty}), E[p])=0$ and that
one has a short exact sequence
\[0 \lra S^{str}(E[p]/F_{\infty}) \lra H^1(G_S(F_{\infty}), E[p]) \lra \bigoplus_{w\in S(F_{\infty})} H^1(F_{\infty,w}, D_w[p])\lra 0. \]
 A standard $\mathbb{F}_p\ps{\Ga}$-corank calculation will show that both $H^1(G_S(F_{\infty}), E[p])$ and $\bigoplus_{w\in S(F_{\infty})} H^1(F_{\infty,w}, D_w[p])$
 have $\mathbb{F}_p\ps{G}$-corank $[F:\Q]$ (see \cite[Theorem 7.1]{HV} and \cite[Theorem 4.1]{OcV03}). It then follows
 that $S^{str}(E[p]/F_{\infty})^{\vee}$ has zero $\mathbb{F}_p\ps{G}$-rank.
  From this, it follows that $X(E/F_{\infty})$ is finitely generated over $\Z_p\ps{H}$.
  \epf


\footnotesize
\begin{thebibliography}{99}
\bibitem{AB} K. Ardakov and K. A. Brown, Primeness, semiprimeness and localisation in Iwasawa algebras,
\textit{Trans. Amer. Math. Soc.} 359(4) (2007) 1499-1515.

\bibitem{A} C. S. Aribam, On the $\mu$-invariant of fine Selmer
groups, \textit{J. Number Theory} 135 (2014) 284-300.

\bibitem{Bh} A. Bhave,
Analogue of Kida's formula for certain strongly admissible
extensions, \textit{J. Number Theory} 122 (2007) 100-120.

\bibitem{C99} J. Coates, Fragments of the $\mathrm{GL}_2$ Iwasawa theory of elliptic curves without complex
multiplication, in \textit{Arithmetic Theory of Elliptic Curves}, ed. C. Viola, Lecture Notes in Math. 1716 (Springer, Berlin, 1999), pp. 1-50.

\bibitem{CG} J. Coates and R. Greenberg, Kummer theory for abelian varieties over
local fields, \textit{Invent. Math.} 124 (1996) 129-174.

\bibitem{CS} J. Coates and R. Sujatha, Fine Selmer groups of elliptic
curves over $p$-adic Lie extensions, \textit{Math. Ann.}
331(4) (2005) 809-839.

\bibitem{CS10} \underline{\hspace{.5in}} , \textit{Galois Cohomology of
Elliptic Curves}, 2nd Ed., Tata Institute of Fundamental Research Lectures on
Mathematics, 88. Published by Narosa Publishing House, New
Delhi; for the Tata Institute of Fundamental Research, Mumbai, 2010.

\bibitem{DSMS} J. Dixon, M. P. F. Du Sautoy, A. Mann and D. Segal, \textit{Analytic Pro-p Groups},
2nd edn, Cambridge Stud. Adv. Math. 38, Cambridge Univ. Press,
Cambridge, UK, 1999.

\bibitem{GW} K. R. Goodearl and R. B. Warfield, \textit{An
introduction to non-commutative Noetherian rings}, London Math. Soc.
Stud. Texts 61, Cambridge University Press, 2004.

\bibitem{G89} R. Greenberg, Iwasawa theory for $p$-adic representations, in
\textit{Algebraic Number Theory--in honor of K. Iwasawa}, ed. J.
Coates, R. Greenberg, B. Mazur and I. Satake, Adv. Std. in Pure
Math. 17, 1989, pp. 97-137.

\bibitem{Gr} \underline{\hspace{.5in}} , Iwasawa theory$-$past and present, in: \textit{Class field Theory$-$its centenary
and prospect}, Adv. Std. in Pure Math. 30, 2001, pp. 335-385.

\bibitem{GV} R. Greenberg and V. Vatsal, On the Iwasawa invariants of elliptic
curves, \textit{Invent. Math.} 142 (2000) 17-63.

\bibitem{HS} Y. Hachimori and R. Sharifi, On the failure of
pseudo-nullity of Iwasawa modules, \textit{J. Alg. Geom.}
14(3) (2005) 567-591.

\bibitem{HV} Y. Hachimori and O.\ Venjakob, Completely faithful
Selmer groups over Kummer extensions. Kazuya Kato's fiftieth
birthday. \textit{Doc. Math.} 2003, Extra Vol., 443-478.

\bibitem{Ho} S. Howson, Euler characteristic as invariants of
Iwasawa modules, \emph{Proc. London Math. Soc.} 85(3) (2002)
634-658.

\bibitem{Ho2} \underline{\hspace{.5in}} , Structure of central torsion Iwasawa
modules, \emph{Bull. Soc. Math. France} 130(4) (2002) 507-535.

\bibitem{Iw} K. Iwasawa, On the $\mu$-invariants of
$\Z_l$-extensions, in: \textit{Number Theory, Algebraic Geometry and
Commutative Algebra, in honour of Yasuo Akizuki}, Kinokuniya, Tokyo,
1973, 1-11.

\bibitem{Iw2} \underline{\hspace{.5in}} , On $\Z_l$-extensions of algebraic
number fields, \textit{Ann.\ of Math.} 98 (1973) 246-326.

\bibitem{Jh} S. Jha, Fine Selmer group of Hida deformations over
non-commutative $p$-adic Lie extensions, \textit{Asian J. Math.}
16(2) (2012) 353-366.

\bibitem{JhS} S. Jha and R. Sujatha, On the Hida deformations
of fine Selmer groups, \textit{J. Algebra} 338 (2011)
180-196.

\bibitem{Lam} T. Y. Lam, \textit{Lectures on Modules and Rings},
Grad. Texts in Math. 189, Springer 1999.

\bibitem{LimMHG}  M. F. Lim, A remark on the $\mathfrak{M}_H(G)$-conjecture
and Akashi series, \textit{Int. J. Number Theory} 11(1)
(2015) 269-297.

\bibitem{LimPS}  \underline{\hspace{.5in}} , On the pseudo-nullity of the dual fine Selmer
   groups, \textit{Int. J. Number Theory} 11(7) (2015) 2055-2063.

\bibitem{LimFine} \underline{\hspace{.5in}} , Notes on the fine Selmer groups,
accepted for publication in \textit{Asian J. Math}.

\bibitem{LimCMu} \underline{\hspace{.5in}} , Comparing the $\pi$-primary
   submodules of the dual Selmer groups,
   accepted for publication in \textit{Asian J. Math}.

\bibitem{LimMurty} M. F. Lim and V. K. Murty, The growth of fine Selmer groups, \emph{J. Ramanujan Math. Soc.}
31(1) (2016) 79-94.

\bibitem{Mazur} B. Mazur, Rational points of abelian varieties in towers of number fields, {\em Invent. Math.} 18 (1972) 183-266.

\bibitem{NSW} J. Neukirch, A. Schmidt and K. Wingberg,
\textit{Cohomology of Number Fields}, 2nd Ed., Grundlehren Math.
Wiss. 323, Springer 2008

\bibitem{Neu} A. Neumann, Completed group algebras without zero divisors,
\textit{Arch. Math.} 51(6) (1988) 496-499.

\bibitem{Oc} Y. Ochi, A remark on the pseudo-nullity conjecture
for fine Selmer groups of elliptic curves, \textit{Comment. Math.
Univ. St. Pauli} 58(1) (2009) 1-7.

\bibitem{OcV02} Y. Ochi and O. Venjakob, On the structure of
Selmer groups over $p$-adic Lie extensions, \textit{J. Alg. Geom.}
11(3) (2002) 547-580.

\bibitem{OcV03}  \underline{\hspace{.5in}} , On the ranks of Iwasawa modules over p-adic Lie extensions, \textit{Math. Proc.
Camb. Phil. Soc.} 135 (2003) 25-43.

\bibitem{RZ} L. Ribes and P. Zalesskii, \textit{Profinite Groups}, Second edition, Ergeb.
Math. Grenzgeb. 40, Springer, 2010.

\bibitem{Sh1}  R. Sharifi, Massey products and ideal class
groups, \textit{J. reine angew. Math.} 603 (2007) 1-33.

\bibitem{Sh2}  \underline{\hspace{.5in}} , On Galois groups of unramified pro-$p$
extensions, \textit{Math. Ann.} 342 (2008) 297-308.

\bibitem{Su} R. Sujatha, Elliptic curves and Iwasawa's
$\mu =0$ conjecture, in: \textit{Quadratic forms, linear algebraic
groups, and cohomology} 125-135, Dev. Math., 18, Springer, New
York, 2010.

\bibitem{V02} O. Venjakob, On the structure theory of the Iwasawa algebra
of a $p$-adic Lie group, \textit{J. Eur. Math. Soc.} 4(3)
(2002) 271-311.

\bibitem{V03} \underline{\hspace{.5in}} , A non-commutative Weierstrass preparation theorem
and applications to Iwasawa theory, \textit{J. reine angew. Math.}
559 (2003) 153-191.

\bibitem{Wu} C. Wuthrich, Iwasawa theory of the fine Selmer group,
\textit{J. Alg. Geom.} 16 (2007) 83-108.

\end{thebibliography}
\end{document}